%% file: main2.tex
\documentclass[runningheads]{llncs}
\usepackage{graphicx}
\usepackage{amssymb}
\usepackage{amsfonts}
\usepackage{amsmath}
\usepackage[english]{babel}

\RequirePackage[l2tabu, orthodox]{nag}
\usepackage[all,warning]{onlyamsmath}

\usepackage{hyperref}
\usepackage{cleveref}
\usepackage{color}
\usepackage{mathtools}
\usepackage{xspace}
\usepackage{fancyvrb}
\usepackage{wrapfig}

\input{macros}

\begin{document}
\title{What's Decidable about Discrete\\ Linear Dynamical Systems?}
%
%
\author{Toghrul Karimov\inst{1} \and Edon Kelmendi\inst{1} \and Jo\"el
  Ouaknine\inst{1}\thanks{Also affiliated with Keble College,
    Oxford as \href{http://emmy.network/}{\texttt{emmy.network}}
    Fellow, and supported by DFG grant 389792660 as part of TRR 248 (see
\url{https://perspicuous-computing.science}).}
  \and James Worrell\inst{2}}
\authorrunning{T. Karimov et al.}
%
\institute{Max Planck Institute for Software Systems, Saarland
  Informatics Campus, Germany \and Department of Computer Science, Oxford University, UK}
\maketitle              
\begin{abstract}
  We survey the state of the art on the algorithmic analysis of
  discrete linear dynamical systems, focussing in particular on
  reachability, model-checking, and invariant-generation questions,
  both unconditionally as well as relative to oracles for the Skolem
  Problem.
  
\end{abstract}

\section{Introduction}

Dynamical systems are a fundamental modelling paradigm in many
branches of science, and have been the subject of extensive research
for many decades. A \emph{(rational) discrete linear dynamical system (LDS)}
in ambient space $\mathbb{R}^d$ is given by a square $d \times d$ matrix $M$ with rational
entries, together with a starting point $x \in \Q^d$.\footnote{All of
  the results we present in this paper carry over to the more
  general setting of \emph{real-algebraic} LDS, whose entries are
  allowed to be real algebraic numbers. Nevertheless, we stick here to
  rationals for simplicity of exposition.}
The
\textit{orbit} of $(M,x)$ is the infinite trajectory
$\mathcal O(M,x) := \langle x, Mx, M^2x,\dots\rangle$. An example of a four-dimensional
LDS is given in Figure~\ref{LDS-example}. Our main focus in the
present paper is on delineating the class of assertions on the orbits of LDS that can be
algorithmically decided.

\begin{figure}
  \label{LDS-example}
\[
M \defeq \begin{pmatrix}
3 & 2&  0 & -5
\\ 
0 &1 &0 &3
\\
0 &4 &3 &13 
\\ 3 &11 &6 &24
\end{pmatrix}
\quad\quad\quad
x \defeq \begin{pmatrix}
1 \\ -1 \\ 2 \\0
\end{pmatrix}
\]
\caption{A four-dimensional discrete linear dynamical system.}
\end{figure}

One of the most natural and fundamental computational questions
concerning linear dynamical systems is the \emph{Point-to-Point
  Reachability Problem}, also known as the \emph{Kannan-Lipton Orbit
  Problem}: given a $d$-dimensional LDS $(M,x)$ together with a point
target $y \in \mathbb{Q}^d$, does the orbit of the LDS ever hit the
target? The decidability of this question was settled affirmatively
in the 1980s in the seminal work of Kannan and
Lipton~\cite{DBLP:conf/stoc/KannanL80,DBLP:journals/jacm/KannanL86}. In
fact, Kannan and Lipton showed that this problem is solvable in
polynomial time, answering an earlier open problem of Harrison from
the 1960s on reachability for linear sequential machines~\cite{Har69}.

Interestingly, one of Kannan and Lipton's motivations was to propose
a line of attack to the well-known \emph{Skolem Problem}, which had itself
been famously open since the 1930s. The Skolem Problem remains unsolved to this day,
although substantial advances have recently been made---more on this shortly.
Phrased in the language of linear dynamical systems, the Skolem
Problem asks whether it is decidable, given $(M,x)$ as above, together
with a $(d-1)$-dimensional subspace $H$ of $\R^d$, to determine if the
orbit of $(M,x)$ ever hits $H$. Kannan and Lipton
suggested that, in ambient space $\R^d$ of arbitrary dimension,
the problem of hitting a low-dimensional subspace might be
decidable. Indeed, this was eventually substantiated by Chonev
\emph{et al.} for linear subspaces of dimension at most
$3$~\cite{DBLP:conf/stoc/ChonevOW13,subspacehit}. 

Subsequent research focussed on the decidability of hitting targets of
increasing complexity, such as
half-spaces~\cite{DBLP:journals/dam/HalavaHH06,DBLP:journals/dam/LaohakosolT09,DBLP:conf/soda/OuaknineW14,OuaknineW14,DBLP:conf/icalp/OuaknineW14a},
polytopes~\cite{DBLP:conf/csr/TarasovV11,polyhit,DBLP:conf/icalp/AlmagorOW17},
and semialgebraic sets~\cite{AlmagorO019,AOW20}. It is also worth
noting that discrete linear dynamical systems can equivalently be
viewed as linear (or affine) simple, branching-free while loops, where
reachability corresponds to loop termination. There is a voluminous
literature on the topic, albeit largely focussing on heuristics and
semi-algorithms (via spectral methods or the synthesis of ranking
functions, in particular), rather than exact decidability results. Relevant papers
include
\cite{DBLP:conf/cav/Tiwari04,DBLP:conf/cav/Braverman06,DBLP:conf/popl/Ben-AmramG13,DBLP:journals/jacm/Ben-AmramG14,DBLP:conf/concur/BradleyMS05,DBLP:journals/sttt/ChenFM15,DBLP:conf/tacas/ColonS01,DBLP:conf/vmcai/PodelskiR04,DBLP:conf/lics/PodelskiR04,DBLP:conf/icalp/HosseiniO019,DBLP:conf/cav/Ben-AmramG17,DBLP:conf/sas/Ben-AmramDG19}. Several
of these approaches have moreover been implemented in software verification
tools, such as Microsoft’s
Terminator~\cite{DBLP:conf/pldi/CookPR06,DBLP:conf/cav/CookPR06}.

In recent years, motivated in part by verification questions for
stochastic systems and linear loops, researchers have begun
investigating more sophisticated decision problems than mere
reachability: for example, the
paper~\cite{DBLP:journals/jacm/AgrawalAGT15} studies approximate LTL
model checking of Markov chains (which themselves can be viewed as
particular kinds of linear dynamical systems),
whereas~\cite{karimov:mfcs} focusses on LTL model checking of
low-dimensional linear dynamical systems with semialgebraic
predicates.\footnote{Semialgebraic predicates are Boolean combinations
  of polynomial equalities and inequalities.} In~\cite{AlmagorKKO021}, the authors solve the
semialgebraic model-checking problem for diagonalisable linear dynamical systems in
arbitrary dimension against prefix-independent MSO\footnote{Monadic
  Second-Order Logic (MSO) is a highly expressive specification
  formalism that subsumes the vast majority of temporal logics
  employed in the field of automated verification, such as Linear
  Temporal Logic (LTL)\@. ``Prefix independence'' is a quality of
  properties that are \emph{asymptotic} in nature---we provide a precise
  definition shortly.} properties, whereas~\cite{POPL22}
investigates semialgebraic MSO model checking of linear dynamical
systems in which the dimensions of predicates are constrained. To
illustrate this last approach, recall the dynamical system $(M,x)$
from Figure~\ref{LDS-example}, and consider the following three
semialgebraic predicates:
\begin{align*}
& P_1(x_1,x_2,x_3,x_4)\defeq x_1 + x_2 + x_3 - x_4  = 0 \land (x_1^3 = x_2^2 \lor x_4 \geq 3x_1^2 +x_2)\\
& P_2(x_1,x_2,x_3,x_4)\defeq x_1 + x_2 + 2x_3-2x_4 = 0 \land x_1^3+x_3^2+x_3 > x_4 \\
& P_3(x_1,x_2,x_3,x_4)\defeq  x_1^4-x_2^2 = 3 \land 2x_3^2 = x_4 \land
                                                                                        x_1^2
                                                                                        -2x_2^3=
                                                                                        4x_3
                                                                                        \, .
\end{align*}

Recall that the ambient space is $\mathbb{R}^4$. We identify the above
predicates with the corresponding subsets of $\mathbb{R}^4$, and wish to
express assertions about the orbit of $(M,x)$ as it traces a
trajectory through $\mathbb{R}^4$. For example (in LTL notation),
\[ \mathbf{G} ( P_1 \Rightarrow \mathbf{F} \lnot P_2) \land \mathbf{F}
  (P_3 \lor \lnot P_1) \] asserts that whenever the orbit visits
$P_1$, then it must eventually subsequently visit the complement of
$P_2$, and moreover that the orbit will eventually either visit $P_3$
or the complement of $P_1$.  The reader will probably agree that
whether or not the above assertion holds for our LDS $(M,x)$ is not
immediately obvious to determine (even, arguably, in
principle). Nevertheless, this example falls within the scope
of~\cite{POPL22}, as the semialgebraic predicates $P_1$, $P_2$, and
$P_3$ are each either contained in
some three-dimensional subspace (this is the case for $P_1$ and
$P_2$), or have intrinsic dimension at most $1$ (this is the case of
$P_3$, which is `string-like', or a curve, as a subset of
$\mathbb{R}^4$). Naturally, we shall return to these notions in due course,
and articulate the relevant results in full details.

A recent and closely related line of inquiry concerns the study of
algebraic model checking of linear dynamical systems~\cite{LOW22}. The
setting is similar to the above, the only difference being that the
allowable predicates are the \emph{constructible} ones, i.e., built
from \emph{algebraic} sets\footnote{Algebraic sets correspond
     to positive Boolean combinations of polynomial equalities.} using
 arbitrary Boolean operations (including
 complementation). The paper~\cite{LOW22} introduces in addition the
 key notion of \emph{Skolem oracle}, which we discuss next.

 \subsection{Skolem Oracles}

 There is an intimate connection between linear dynamical systems and
linear recurrence sequences. An \emph{(integer) linear recurrence sequence (LRS)}
$\boldsymbol{u} = \langle u_n \rangle_{n=0}^\infty$ is an infinite sequence of
integers satisfying
\begin{equation}
\label{eq:RECUR}
u_{n+d} = c_1u_{n+d-1} + \cdots + c_{d-1}u_{n+1}+ c_du_{n} 
\end{equation}
for all $n \in \mathbb{N}$, where the coefficients $c_1,\ldots,c_d$
are integers and $c_d \neq 0$.  We say that the above
recurrence has \emph{order} $d$. We moreover say that an LRS is
\emph{simple} if the characteristic polynomial\footnote{The
  characteristic polynomial associated with recurrence~(\ref{eq:RECUR})
  is $X^d-c_1X^{d-1}- \ldots -c_d$.}  of its minimal-order recurrence
has no repeated roots. The sequence of Fibonacci numbers
$\langle f_n \rangle_{n=0}^\infty = \langle 0,1,1,2,3,5,\ldots
\rangle$, which obeys the recurrence $f_{n+2}=f_{n+1}+f_n$, is perhaps
the most emblematic LRS, and also happens to be simple. 

The celebrated theorem of Skolem, Mahler, and Lech (see~\cite{Everest/0035746})
describes the structure of the set $\{ n \in \mathbb{N} : u_n = 0\}$
of zero terms of an LRS as follows:
\begin{theorem}
  \label{thm:SML}
  Given a linear recurrence sequence
  $\boldsymbol{u}=\langle u_n\rangle_{n=0}^\infty$, its set of zero
  terms is a semilinear set, i.e., it consists of a union of finitely
  many full arithmetic progressions,\footnote{A full arithmetic
    progression is a set of non-negative integers of the form
    $\{a + bm : m \in \N\}$, with $a,b \in \N$ and $a < b$.}  together with a finite set.
\end{theorem}

As shown by Berstel and Mignotte~\cite{berstel76_deux_des_suites}, in
the above one can effectively extract all of the arithmetic
progressions; we refer herein to the corresponding procedure as the
`Berstel-Mignotte algorithm'. Nevertheless, how to compute the leftover finite set of
zeros remains open, and is easily seen to be equivalent to the \emph{Skolem
Problem}: given an LRS $\boldsymbol{u}$, does $\boldsymbol{u}$ contain
a zero term?

The paper~\cite{LOW22} therefore introduces the notion of a \emph{Skolem oracle}: given
an LRS $\boldsymbol{u}=\langle u_n\rangle_{n=0}^\infty$, such an
oracle returns the finite set of indices of zeros of $\boldsymbol{u}$
that do not already belong to some infinite arithmetic progression of
zeros. Likewise, a \emph{Simple-Skolem oracle} is a Skolem oracle
restricted to simple LRS\@.

As mentioned earlier, the decidability of the Skolem Problem is a
longstanding open
question~\cite{Everest/0035746,DBLP:journals/siglog/OuaknineW15}, with
a positive answer for LRS of order at most $4$ known since the
mid-1980s~\cite{Tijdeman1984,vereshchagin1985problem}.  Very recently,
two major conditional advances on the Skolem Problem have been made,
achieving decidability subject to certain classical number-theoretic
conjectures: in~\cite{LLN22}, Lipton \emph{et al.}\ established
decidability for LRS of order $5$ assuming the \emph{Skolem
  Conjecture} (also known as the \emph{Exponential Local-Global
  Principle}); and in~\cite{BLN22}, Bilu \emph{et al.} showed
decidability for simple LRS of arbitrary order, subject to both the
Skolem Conjecture and the \emph{$p$-adic Schanuel Conjecture} (we
refer the reader to~\cite{BLN22} for the precise definitions and details). It is
interesting to note that in both cases, the procedures in question
rely on the conjectures \emph{only} for termination; correctness is
unconditional. In fact, these procedures are \emph{certifying
  algorithms} (in the sense of~\cite{MMN11}) in that, upon
termination, they produce an independent certificate (or witness) that
their output is correct. Such a certificate can be checked
algorithmically by a third party with no reliance on any unproven conjectures. The
authors of~\cite{BLN22} have implemented their algorithm within the
\textsc{skolem} tool, available
online.\footnote{\texttt{https://skolem.mpi-sws.org/} .}

In view of the above, Simple-Skolem oracles \emph{can} be implemented
with unconditional correctness, and guaranteed termination subject to
the Skolem and $p$-adic Schanuel conjectures. Whether full Skolem oracles
can be devised is the subject of active research (see, e.g.,~\cite{LOW21,LOW22a}); at the time of
writing, to the best of our knowledge, no putative procedure is even
conjectured in the general (non-simple) case.

\subsection{Paper Outline}

Questions of reachability and model checking for linear dynamical
systems constitute one of the central foci of this paper. In
Section~\ref{sec:model-checking}, we cover the state of the art and
beyond, both unconditionally and relative to Skolem oracles. We paint
what is essentially a complete picture of the landscape, in each
situation establishing either decidability (possibly conditional on
Skolem oracles), or hardness with respect to longstanding open
problems. An important theme in the classical theory of dynamical
systems concerns the study of \emph{asymptotic} properties (e.g., stability,
convergence, or divergence of orbits), and we therefore consider
both MSO along with its \emph{prefix-independent} fragment piMSO\@.
Section~\ref{sec:pseudo} then focusses on
questions of robustness through the notion of pseudo-orbit. In
Section~\ref{sec:invariants}, we discuss the algorithmic synthesis of
inductive invariants for linear dynamical systems, and
Section~\ref{sec:initial-sets} examines the situation in which orbits
originate from an initial set rather than a single point. Finally, 
Section~\ref{sec:future} concludes with a brief summary and a glimpse
of several research directions.

\section{Model Checking}
\label{sec:model-checking}

Throughout this section, we assume familiarity with the rudiments of Monadic Second-Order
Logic (MSO); an excellent reference is the text~\cite{BGG97}.

Let us work in fixed ambient space $\R^d$, and consider a
$d$-dimensional LDS $(M,x)$ (i.e., $M \in \Q^{d\times d}$ and $x \in \Q^d$).
Recall that the orbit $\mathcal{O} = \mathcal{O}(M,x)$ of our LDS is the infinite
sequence $\langle x, Mx, M^2x, \ldots \rangle$ in $\Q^d$. Let us write
$\mathcal{O}[n]$ for the $n$th term of the orbit.

Given an MSO formula $\varphi$ over the collection of semialgebraic
predicates $\PR= \{P_1, \ldots, P_m\}$, where each
$P_i \subseteq \R^d$, the model-checking problem consists in
determining whether the orbit (more precisely, the characteristic word
$\alpha \in \left(2^\PR\right)^\omega$ of the orbit
$\mathcal{O}(M, x)$ with respect to $\PR$, where $P_i \in \alpha[n]$
iff $\mathcal{O}[n] \in P_i$) satisfies
$\varphi$. Reachability problems for LDS constitute special cases of
the model-checking problem, and already the questions of determining
whether a given orbit reaches a hyperplane (Skolem Problem) or a
halfspace (the \emph{Positivity Problem}
\cite{OuaknineW14}) are
longstanding open problems in number theory, couched in the language of
linear dynamical systems. Recent research has, however, succeeded in
uncovering several important decidable subclasses of the
model-checking problem and in demarcating the boundary between what is
decidable and what is hard with respect to longstanding open mathematical problems.

In order to present the main results, we require some further
definitions to specify the classes of predicates that are allowed within MSO formulas.
Let us write $\SA \subseteq 2^{\R^d}$ and $\CO \subseteq
2^{\R^d}$ to denote respectively the
collections of all semialgebraic subsets of $\R^d$ and of all
constructible subsets of $\R^d$.\footnote{Recall that $\CO$ is the
smallest set containing all algebraic subsets of $\R^d$, and which is
closed under finite union, finite intersection, and complement. (The
terminology of ``constructible'' originates from algebraic geometry.)}
We also define the collection $\T \subseteq 2^{\R^d}$ of \emph{tame}
sets as follows: $\T$ comprises all semialgebraic subsets of
$\R^d$ that are either contained in a three-dimensional
subspace of $\R^d$, or that have intrinsic dimension at most
one.\footnote{The intrinsic dimension of a semialgebraic set is
  formally defined via cell decomposition; intuitively,
  one-dimensional semialgebraic sets can be viewed as ‘strings’ or
  ‘curves’, whereas zero-dimensional semialgebraic sets are finite
  collections of singleton points.}  Moreover, $\T$ is defined to be the smallest
such set which is in addition closed under finite union, finite
intersection, and complement. Finally,
we define the set $\T \oplus \CO$ to be the smallest superset of
$\T \cup \CO$ that is closed under finite union, finite
intersection, and complement.

Note that all of $\SA$, $\CO$, $\T$, and $\T \oplus \CO$ are closed
under all Boolean operations; this is in keeping with their intended
use as collections of predicates for MSO formulas, bearing in mind
that MSO itself possesses all Boolean operators.

The motivation for considering the collection $\T$ of tame predicates
 has origins in the results of
\cite{DBLP:conf/stoc/ChonevOW13,karimov:mfcs,AlmagorO019,BaierFJLLOPWW2021Orbit}. A
common theme is that for tame predicates, the proofs that establish
how to decide reachability also provide one with a means of
representing, in a finitary manner, all the time steps at which the
orbit of a given LDS is in a particular predicate set $T$. The authors of \cite{POPL22}
show how to combine these representations (one for each predicate) to
obtain structural information about the characteristic word $\alpha$
that is sufficient for determining whether a deterministic automaton
$\aut$ accepts $\alpha$, leading to the following.

\begin{theorem}
  \label{thm:mc1}
  Let $(M,x)$ be an LDS, $\PR = \{T_1, \ldots, T_m\} \subseteq \T$ be a collection of
  tame predicates and $\varphi$ be an MSO formula over
  $\PR$. It is decidable whether $(M,x) \vDash \varphi$ (i.e., whether
  the characteristic word $\alpha$ of the orbit $\mathcal{O}(M,x)$ with
  respect to $\PR$ satisfies $\varphi$).
\end{theorem}

Let us note in passing that Theorem~\ref{thm:mc1} subsumes the
decidability of the Kannan-Lipton Point-to-Point Reachability Problem,
since points are singleton sets and the latter are evidently tame.

It is also worth pointing out, absent other restrictions, that this delineation of the decidable fragment of
the model-checking problem is tight as trying to expand the definition
of tame predicates runs
into open problems already for formulas that describe mere
reachability properties. In particular, the Skolem Problem in
dimension 5 is open and can be encoded (i)~as a reachability problem
with a four-dimensional LDS and a three-dimensional affine subspace
\cite{DBLP:conf/stoc/ChonevOW13} (that is, in general, not contained
in a three-dimensional linear subspace) and (ii)~as a reachability
problem with a target of intrinsic dimension two
\cite{BaierFJLLOPWW2021Orbit}.

To sidestep these obstacles, in \cite{AlmagorKKO021} the authors
restrict $\varphi$ to formulas that define \emph{prefix-independent}
properties. A property is
prefix-independent if the infinite words that satisfy it are closed
under the operations of insertion and deletion of finitely many
letters.\footnote{It is interesting to note that whether an MSO
  formula $\varphi$ is prefix-independent or not is decidable.
To see this, for $A=(Q,q_0,\Sigma,\Delta,F)$ a deterministic M\"uller
automaton, define $A(q)$, for $q\in Q$, to be the same as $A$, except
that the initial state of $A(q)$ is $q$ (rather than $q_0$). We say
that a deterministic M\"uller automaton $A$ (as above) is 
\emph{prefix-independent} if, for all $q\in Q$ that are reachable from $q_0$, $A(q)$
recognises the same language as $A$. Write $L(A)$ to denote the
language recognised by $A$. It is now straightforward to show
that $A$ is prefix-independent iff $L(A)$ is prefix-independent. Since
any MSO formula is encodable as a deterministic M\"uller automaton,
and equality of $\omega$-regular languages is decidable,
the desired decidability result follows.} Such properties capture behaviours that are intrinsically
asymptotic in nature (for example: ``does the orbit enter $P$
infinitely often?''); note that the property of
reachability is \emph{not} prefix-independent. The main theorem of
\cite{AlmagorKKO021} in this direction concerns \emph{diagonalisable}
linear dynamical systems:\footnote{An LDS $(M,x)$ is
  \emph{diagonalisable} if the matrix $M$ is diagonalisable (over
  $\C$). In a measure-theoretic sense, most LDS are diagonalisable.}

\begin{theorem}
  \label{thm:mc2}
  Let $(M, x)$ be a diagonalisable LDS and $\varphi$ be a
  prefix-independent MSO formula over a collection of semialgebraic
  predicates $S_1,\ldots,S_m \in \SA$. It is decidable whether $(M,x)
  \vDash \varphi$.
      \end{theorem}
      
      Note in the above that the semialgebraic predicates are entirely
      unrestricted (in particular, not required to be tame).
      However, the restrictions to prefix-independent formulas and
      to diagonalisable systems both again turn out to be essential. Since
      the Skolem Problem is open for diagonalisable systems (in
      dimensions $d \geq 5$), the (non-prefix-independent) model-checking problem for
      diagonalisable LDS is Skolem-hard already for four-dimensional
      systems and affine subspace targets, as discussed earlier. On
      the other hand, if we allow non-diagonalisable systems, then the
      problem of determining whether the orbit of an LDS is eventually
      trapped in a given half-space $H$ (known as the
      \emph{Ultimate Positivity Problem}, corresponding to the
      prefix-independent formula
      $\varphi = \mathbf{F} \, \mathbf{G} \, H$) is hard with respect
      to certain longstanding open problems in Diophantine
      approximation
      \cite{OuaknineW14}.

     This last observation however suggests that it might be possible
     to orchestrate a trade-off between the type of LDS under
     consideration and the class of allowable specification
     predicates. Indeed, it turns out that one can lift the
     restriction to diagonalisable LDS if one agrees to restrict the class
     of predicates:

     \begin{theorem}
 \label{thm:mc3}
  Let $(M, x)$ be an LDS and $\varphi$ be a
  prefix-independent MSO formula over a collection of 
  predicates $P_1,\ldots,P_m \in \T \oplus \CO$. It is decidable whether $(M,x)
  \vDash \varphi$.
\end{theorem}

Theorem~\ref{thm:mc3} goes beyond both Theorem~\ref{thm:mc1} (in that
constructible predicates are allowed into the mix) as well
as~\cite[Thm.~7.3]{LOW22} (in that tame predicates are allowed). 

Let us provide a brief proof sketch of Theorem~\ref{thm:mc3}. Let $(M,x)$
and $\varphi$ be as above. Observe first
(as a straightforward exercise) that any $P \in \T \oplus \CO$ can be
written in conjunctive normal form,
i.e., as an expression of the form
$P = \bigcap_{i=1}^a \bigcup_{j=1}^b B_{i,j}$, where each $B_{i,j}$ is
either a tame predicate or a constructible predicate. Without loss of
generality, one may therefore assume that each predicate appearing in 
$\varphi$ is either tame or constructible.

We now invoke~\cite[Prop.~5]{LOW22} to conclude that, for each constructible
predicate $B_\ell$ appearing in $\varphi$, the Boolean-valued word
$\alpha_\ell$ tracking the passage of the orbit of $(M,x)$ through
$B_\ell$ is ultimately periodic. Moreover, thanks to the
Berstel-Mignotte algorithm, the attendant arithmetic progressions can
all be effectively elicited. For each such $\alpha_\ell$, one can
therefore construct a (fully) periodic word $\alpha'_\ell$ which differs
from $\alpha_\ell$ in at most finitely many places. In  other words,
for all sufficiently large $n$, $M^nx \in B_\ell$ iff $\alpha'[n] = \mathit{true}$.
Being periodic, $\alpha'_\ell$ can in turn be described by an MSO
subformula $\psi_\ell$.
Let us therefore replace within $\varphi$
every occurrence of a constructible predicate $B_\ell$ by the
subformula $\psi_\ell$, obtaining in this process a
new MSO formula $\varphi'$ that comprises exclusively tame
predicates. As $\varphi$ is prefix-independent, it is immediate that 
$(M,x) \vDash \varphi$ iff $(M,x) \vDash \varphi'$. But the latter is of
course decidable thanks to Theorem~\ref{thm:mc1}, concluding the proof
sketch of Theorem~\ref{thm:mc3}.

Once again, Theorem~\ref{thm:mc3} is tight: in Appendix~\ref{appA}, we
show that the ability to solve the model-checking problem for
prefix-independent MSO specifications making use of semialgebraic
predicates in ambient space $\R^4$ would necessarily entail major
breakthroughs in Diophantine approximation.

Let us now turn to the question of the extent to which
the above results can be enhanced through the use of Skolem oracles. The
key result is as follows, in effect enabling us to drop the
restriction of prefix-independence from Theorem~\ref{thm:mc3}:

 \begin{theorem}
 \label{thm:mc4}
 Let $(M, x)$ be an LDS and $\varphi$ be an MSO formula over a collection of 
  predicates $P_1,\ldots,P_m \in \T \oplus \CO$. It is decidable whether $(M,x)
  \vDash \varphi$, subject to the existence of a Skolem oracle.

  The same result also holds for diagonalisable LDS, assuming the
  existence of a Simple-Skolem oracle.
\end{theorem}

The proof of Theorem~\ref{thm:mc4} is similar to that of
Theorem~\ref{thm:mc3}; we provide a brief sketch below.

Let $(M,x)$ and $\varphi$ be as above, and assume, thanks to the
representation of ($\T \oplus \CO$)-predicates in conjunctive normal
form, that every predicate occurring in $\varphi$ is either tame of
constructible. Thanks to~\cite[Cor.~6]{LOW22}, for each constructible
predicate $B_\ell$ appearing in $\varphi$, the Boolean-valued word
$\alpha_\ell$ tracking the passage of the orbit of $(M,x)$ through
$B_\ell$ is effectively ultimately periodic (this requires the use of
a Skolem oracle, or a Simple-Skolem oracle if $M$ is
diagonalisable). In other words, we have a finitary exact
representation of $\alpha_\ell$, and can therefore describe it via an
MSO subformula $\psi_\ell$.

We can now replace within $\varphi$ every occurrence of a
constructible predicate $B_\ell$ by its corresponding subformula
$\psi_\ell$, obtaining in this process an equivalent MSO formula
$\varphi'$ that comprises exclusively tame predicates. The desired
result then immediately follows from Theorem~\ref{thm:mc1}, concluding
the proof sketch of Theorem~\ref{thm:mc4}.

As noted earlier, Simple-Skolem oracles can be implemented into
provably correct certifying procedures which terminate subject to
classical number-theoretic conjectures~\cite{BLN22}. Let us therefore
separately record an important corollary:

\begin{corollary}
\label{cor:mc4}
Let $(M, x)$ be a diagonalisable LDS and $\varphi$ be an MSO formula
over a collection of predicates $P_1,\ldots,P_m \in \T \oplus \CO$. It
is decidable whether $(M,x) \vDash \varphi$, assuming the Skolem
Conjecture and the p-adic Schanuel Conjecture. Moreover, correctness
of the attendant procedure is unconditional, and independent
correctness certificates can be produced upon termination.
\end{corollary}

Let us point out that Theorem~\ref{thm:mc4} is, once again, tight. For
arbitrary LDS, a similar argument as that put forth in
Appendix~\ref{appA} applies, since it is not known (or even believed) that
Skolem oracles are of any use in tackling Ultimate Positivity
problems. For diagonalisable LDS, we can invoke the order-10
Simple Positivity Problem, which remains open to this day
(see~\cite{OuaknineW14}); it can be modelled straightforwardly as a
half-space reachability problem in ambient space $\R^{10}$ (or even $\R^9$, by
considering an affine half-space). In fact, critical unsolved cases of order-10
 Positivity can even be formulated as semialgebraic reachability
problems in ambient space $\R^4$; we omit the details in the interests of
space and simplicity of exposition.

We summarise the main results of this section in
Figure~\ref{sec_2:summary} below.

  \begin{figure}
  \label{sec_2:summary}
\begin{center}    
  \begin{tabular}{ | r | c | c | }
 \multicolumn{3}{c}{\textbf{Diagonalisable LDS:}}\vspace{0.5ex} \\
    \hline
    & \ unconditional\ \  & \ Simple-Skolem oracle\ \  \\
    \hline
 \ piMSO\ \ & $\SA$ (Thm.~\ref{thm:mc2}) & $\SA$ (Thm.~\ref{thm:mc2}) \\
 \hline
 MSO\ \  & $\T$ (Thm.~\ref{thm:mc1})  & $\T \oplus \CO$ (Thm.~\ref{thm:mc4})\\
  \hline
  \end{tabular}

  \vspace{3ex}

 \begin{tabular}{ | r | c | c | }
 \multicolumn{3}{c}{\textbf{Arbitrary LDS:}}\vspace{0.5ex} \\
    \hline
    & \ unconditional\ \  & \ Skolem oracle\ \  \\
    \hline
 \ piMSO\ \ &\ $\T \oplus \CO$ (Thm.~\ref{thm:mc3})\ \ &\ $\T \oplus \CO$
                                                     (Thm.~\ref{thm:mc3}
                                                     or
                                                         Thm.~\ref{thm:mc4})\ \ \\
 \hline
 MSO\ \  & $\T$ (Thm.~\ref{thm:mc1})  & $\T \oplus \CO$ (Thm.~\ref{thm:mc4})\\
  \hline
  \end{tabular}
  
  \end{center}

\caption{A summary of the decidable model-checking fragments for both diagonalisable
  and arbitrary linear dynamical systems. The prefix-independent
  fragment of MSO is denoted piMSO. $\SA$ is the collection of
  semialgebraic predicates, $\T$ is the collection of tame predicates
  (Boolean closure of semialgebraic sets that either are contained in a
  three-dimensional subspace, or have intrinsic dimension at most
  one), $\CO$ is the collection of constructible predicates (Boolean
  closure of algebraic sets), and $\T \oplus \CO$ is the Boolean
  closure of $\T \cup \CO$. The right-hand columns in both tables
  assume access to Skolem or Simple-Skolem oracles.}
\end{figure}

Taken together, Theorems~\ref{thm:mc1}--\ref{thm:mc4}, along with Corollary~\ref{cor:mc4},
not only subsume---to the best of our knowledge---all existing results
regarding model-checking and reachability problems for discrete linear
dynamical systems, but moreover paint an essentially complete picture
of what is (even in principle) feasible, barring major breakthroughs
in longstanding open problems. It is noteworthy that, in this
characterisation, there appears to be very little difference between
being able to decide mere reachability for a given class of predicates,
and being able to decide the whole of MSO over the same class of predicates.

\section{Pseudo-Reachability and Robustness}
\label{sec:pseudo}

In this section we discuss decision problems about pseudo-orbits that
are related to robustness of computation. Given an LDS $(M,x)$, recall
that the orbit of $x$ under $M$ is the sequence
$\langle x, Mx, M^2x,\dots\rangle$. We say that the sequence
$\langle x_n : n \in \nat \rangle$ is an
\emph{$\varepsilon$-pseudo-orbit} of $x$ under $M$ if $x_0 = x$ and
$x_{n+1} =M x_n + d_n$ for some perturbation $d_n$ with
$||d_n|| < \epsilon$. The \emph{pseudo-orbit} of $x$ under $M$ is then
defined as the set of points that are reachable from $x$ via an
$\varepsilon$-pseudo-orbit for every $\epsilon>0$. This notion of an
($\varepsilon$-)pseudo-orbit, introduced and studied by Anosov~\cite{Anosov},
Bowen~\cite{Bowen1} and Conley~\cite{Conley78}, is an important
conceptual tool in dynamical systems. From the computational
perspective, an $\varepsilon$-pseudo-orbit can be viewed as a
trajectory after a rounding error of magnitude at most $\epsilon$ is
applied at each step.

Given these definitions, we can consider the reachability and
model-checking problems for pseudo-orbits. A natural analogy to the
Kannan-Lipton Orbit Problem is the \emph{Pseudo-Orbit Problem}, which
is to determine whether a target point $y$ belongs to the pseudo-orbit
of $x$ under $M$. In \cite{dcosta_et_al:LIPIcs.MFCS.2021.34} the
authors show that, just like the Orbit Problem, the Pseudo-Orbit
Problem is decidable in polynomial time.  Generalising from points to
sets, let us say that a target set $T$ is pseudo-reachable if for
every $\epsilon > 0$ there exists an $\varepsilon$-pseudo-orbit of $x$
under $M$ that reaches $T$. We can then define the \emph{Pseudo-Skolem
  Problem} and the \emph{Pseudo-Positivity Problem} to be the
pseudo-reachability problems with a hyperplane and a halfspace as
target sets, respectively. Surprisingly,
\cite{dcosta_et_al:LIPIcs.MFCS.2021.34} shows that both of these
problems are in fact decidable! Moreover, \cite{DKMOSW22} establishes
decidability of pseudo-reachability with arbitrary
semialgebraic targets for diagonalisable linear dynamical systems.

Inspired by the above results, one might consider the model-checking
problem for pseudo-orbits, namely the problem of determining, given
$(M, x)$ and a formula $\varphi$, whether for every $\varepsilon > 0$,
there exists an $\varepsilon$-pseudo-orbit that satisfies
$\varphi$. After all, as discussed in the preceding section, for
genuine orbits the fragments of the reachability problem and the full
MSO model-checking problem that are known to be decidable (i.e., the
restrictions on the class of predicates and the property $\varphi$ that make
the problems decidable) are essentially the same. This optimism is,
however, quickly shattered by the following observation. Let $H$ be a
closed halfspace and $\varphi$ be the property $\mathbf{G} \, H$
(``the trajectory always remains inside $H$''). Then the pseudo-orbits
satisfy $\varphi$ (in the sense defined above) if and only the
(genuine) orbit satisfies $\varphi$. The problem of determining
whether the orbit $\mathcal{O}(M,x)$ always remains in $H$, is
however, equivalent to the problem of determining whether the orbit
ever hits an open halfspace, which itself is the Positivity
Problem (a longstanding open question).

\section{Invariant Generation}
\label{sec:invariants}

In the absence of fully general algorithms to decide whether the orbit
of a given LDS reaches targets of arbitrary forms, much effort has
been expended on sound---but possibly incomplete---techniques,
and particularly on constructing certificates of (non-)reachability. This
splits into two broad lines of attack: ranking functions and
invariants. The former are certificates of reachability, demonstrating
that progress is being made towards the target. Inductive invariants
are, on the other hand, certificates of non-reachability, 
establishing that the orbit will not reach the target by enclosing the
former within a set that is itself disjoint from the latter. We focus
in this section on the algorithmic generation of invariants.

More precisely, a set $\mathcal I\subseteq \rel^d$ is said to be an
\emph{inductive invariant} of $(M,x)$ if it contains $x$
($x\in\mathcal I$), and is stable under $M$, that is:
\begin{align*}
  M\mathcal I\defeq \set{My : y\in\mathcal I}\subseteq \mathcal I \, . 
\end{align*}
Clearly there are some trivial invariants, such as $\rel^d$ and the
orbit $\mathcal O(M,x)$ itself. They are not particularly useful in
the sense that the ambient space $\rel^d$ is never disjoint from
whatever target might be under consideration, whereas for various
classes of targets (such as hyperplanes or half-spaces; or more
generally arbitrary semialgebraic sets) we do not in general know how
to decide whether $\mathcal O(M,x)$ is disjoint from the target.
Hence one does not seek \emph{any} invariant, but rather an invariant
that can be algorithmically established to be disjoint from the
target.
\begin{center}
\includegraphics[width=.7\textwidth]{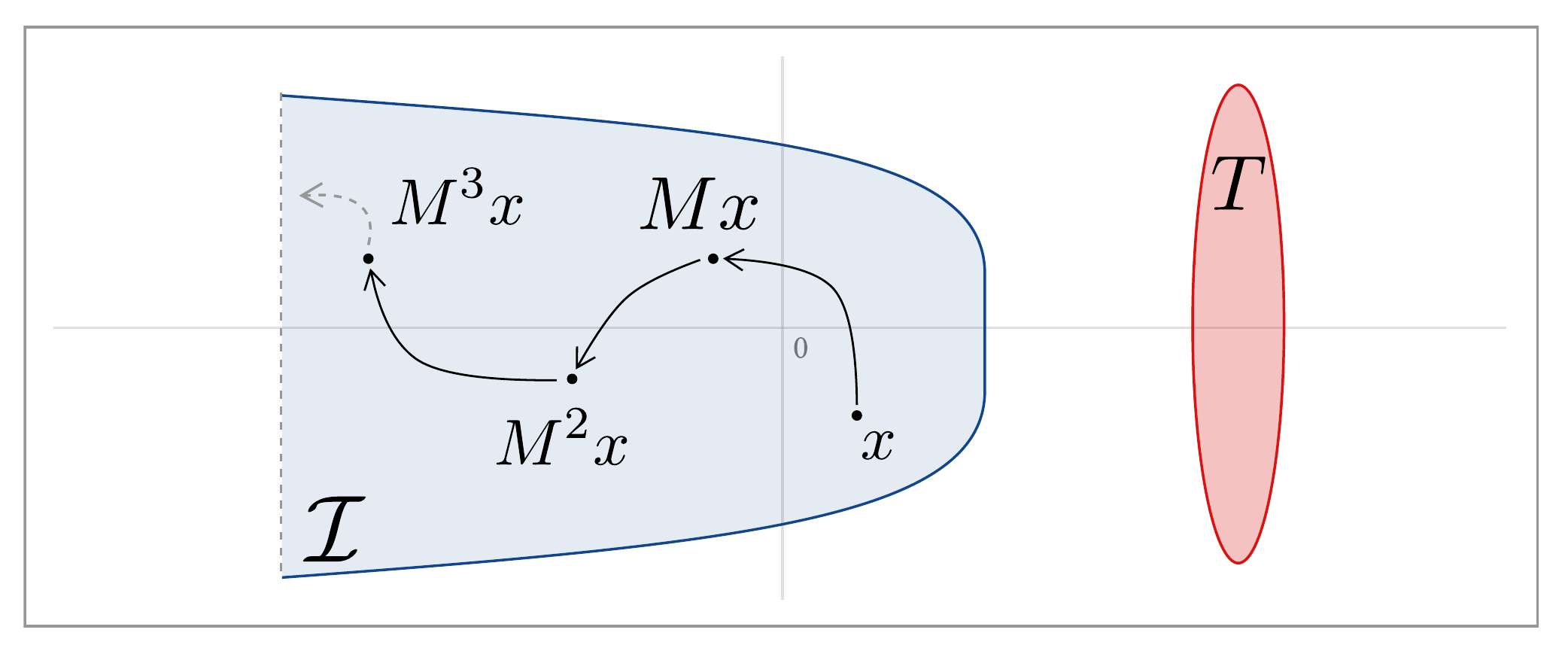}
\end{center}
We therefore seek a sufficiently large, or expressive, class of invariants $\mathcal F$ which moreover exhibits favourable algorithmic properties. A natural family to consider is the collection of all semialgebraic sets. We now have:

\begin{theorem}[{\cite{DBLP:journals/tocl/AlmagorCOW22}}]
  \label{th:invariants}
  Given an LDS $(M,x)$ in ambient space $\rel^d$, together with a semialgebraic target $T\subseteq \rel^d$, it is decidable whether there exists a semialgebraic invariant of $(M,x)$ that is disjoint from $T$. 
\end{theorem}
Furthermore, the algorithm explicitly constructs the invariant when it exists, in the form of a Boolean combination of polynomial inequalities.

Theorem~\ref{th:invariants} holds for an even larger class
$\mathcal F$, namely that of o-minimal sets. We give an informal
definition. Recall the contents of Tarski's quantifier-elimination
theorem, to the effect that semialgebraic subsets of $\rel^d$ are
closed under projections. Moreover, semialgebraic subsets of $\rel$
are quite simple: they are finite unions of intervals. Other families
of sets that enjoy these two properties exist, notably those definable
in the first-order theory of the real numbers augmented with a symbol
for the exponential function, an important result due to
Wilkie~\cite{wilkie1996model}. Structures of $\rel^d$ that are induced
by such logical theories are called \emph{o-minimal}, and an
\emph{o-minimal set} is a set that belongs to such a
structure~\cite{van1998tame}. These include semialgebraic sets, as
well as sets definable in the first-order theory of the reals with
restricted analytic functions.

In \cite{DBLP:conf/icalp/AlmagorCO018}, it is shown that it is decidable whether there exists an o-minimal invariant (for a given LDS $(M,x)$) that is disjoint from a semialgebraic target $T$; and moreover, when such an invariant exists, it is always possible to exhibit one that is in fact semialgebraic~\cite{DBLP:journals/tocl/AlmagorCOW22}. Once again, these results are effective, and the invariants can always be explicitly produced.

\begin{wrapfigure}[13]{r}{.3\textwidth}
  \begin{center}
    \vspace{-1.3cm}
    \includegraphics[width=.3\textwidth]{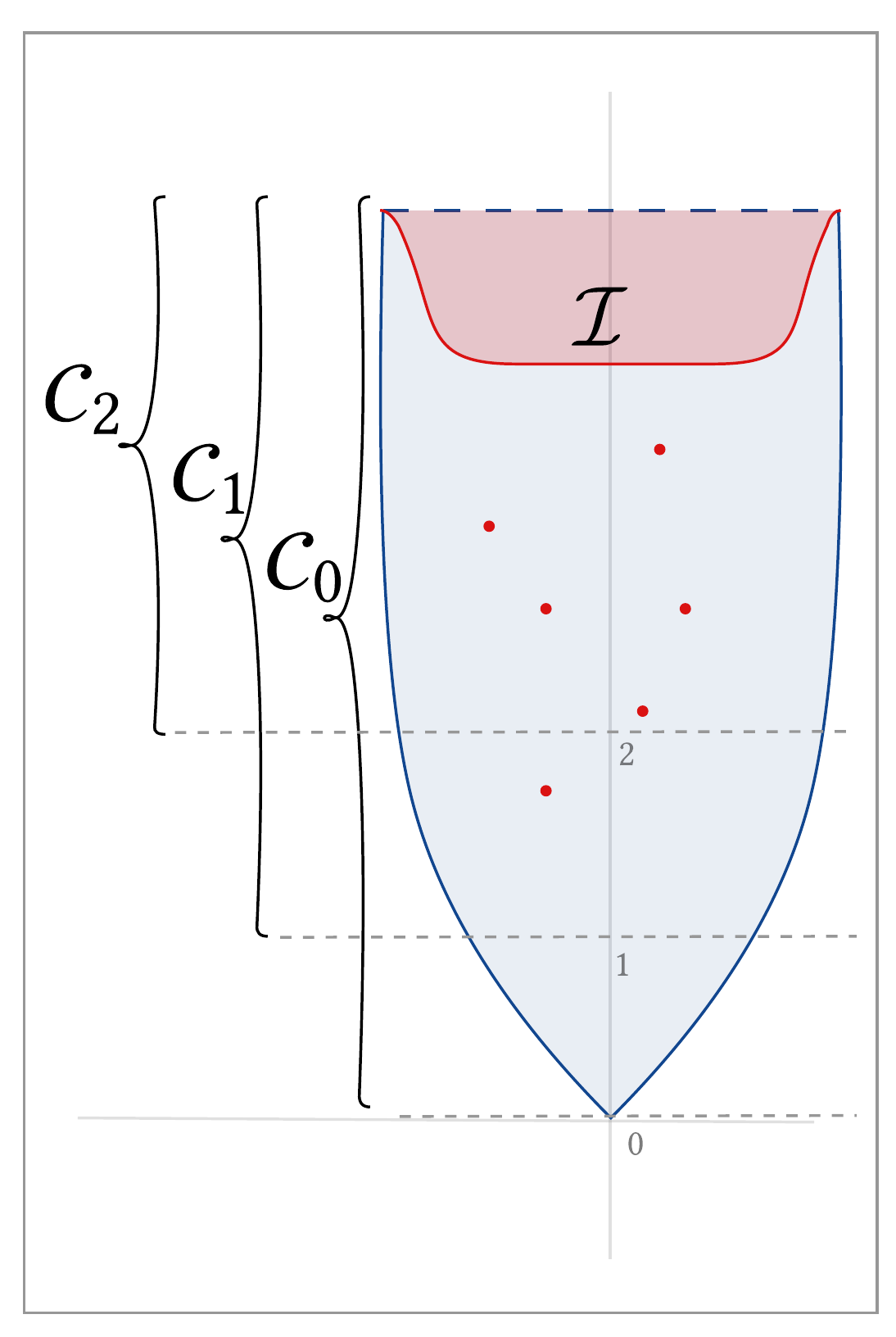}
  \end{center}
\end{wrapfigure}
Varying the class of invariants and the class of targets gives rise to
a number of natural questions that---for the most part---remain
unexplored. Let us however mention a further desirable property
enjoyed both by the class of o-minimal sets and that of semialgebraic
sets: in either case, they admit \emph{minimal families} of
invariants, a notion which we now explain. Let $\mathcal{F}$ be a
class of sets---either the class of o-minimal sets or that of
semialgebraic sets. It is easy to verify that, in general, $\mathcal{F}$ does not
possess minimal invariants.
Nevertheless,
\cite{DBLP:conf/icalp/AlmagorCO018,DBLP:journals/tocl/AlmagorCOW22}
show how to produce a sequence of $(M,x)$-invariants
$\langle \mathcal C_k : k\in\nat \rangle$, all belonging to
$\mathcal{F}$, and such that $\mathcal C_{k+1}\subset \mathcal
C_k$. It can moreover be shown that, given any $(M,x)$-invariant
$\mathcal{I} \in \mathcal{F}$, it is always the case that
$\mathcal{I}$ contains one of the $\mathcal{C}_k$, and \emph{ipso
  facto} also all $\mathcal{C}_j$ for $j \geq k$.


\section{Semialgebraic Initial Sets}
\label{sec:initial-sets}

Up until now we have exclusively considered problems concerning the orbit $\mathcal O(M,x)$ of a \emph{single} initial point $x$. It is natural to ask whether the algorithmic problems which we have discussed remain solvable if one instead considers an entire \emph{set} of initial points $S\subseteq \rel^d$. Unfortunately the answer is negative. 
We sketch below the proof of the undecidability of a natural model-checking problem.

\begin{theorem}
  \label{th:undecidable}
  The following problem is undecidable. Given a natural number $k\in\nat$, a semialgebraic set $S\subseteq \rel^d$, a $d\times d$ rational matrix $M$, and a hyperplane $H$ in $\rel^d$ (having rational normal vector),
determine whether there exists $x\in S$ such that the orbit generated by $(M,x)$ hits $H$ at least $k$ times. 
\end{theorem}
In symbols, whether there is some $x\in S$ such that
\begin{align*}
  \lvert\mathcal O(M,x)\cap H\rvert \ge k
\end{align*}
is undecidable.

It is worth noting that all the problems that we have discussed so far
(including the Skolem and Positivity Problems) are not known to be
undecidable, and are in fact widely conjectured to be decidable. It is
therefore perhaps somewhat surprising that this natural generalisation
of our setting immediately leads to undecidability.

The proof of Theorem~\ref{th:undecidable} proceeds by reduction from a
variant of Hilbert's tenth problem. Recall that Hilbert's tenth
problem asks whether a given polynomial
$P\in\intg[Y_1,Y_2,\ldots,Y_{d-1}]$ with integer coefficients and
$d-1$ variables has a root with all unknowns taking integer
values. This problem is undecidable, as shown by Davis, Putnam, Robinson, and Matiyasevich \cite{hilbert}.

The variant that we will reduce from asks whether the polynomial has roots with the unknowns being distinct natural numbers. It is straightforward to show that this variant is also undecidable.

Let $d\in\nat$, $d>1$, and $P\in\intg[Y_1,Y_2,\ldots,Y_{d-1}]$ be an arbitrary polynomial. We define the subset $S\subseteq \rel^{d}$ via a formula of the first-order theory of the reals:
\begin{align*}
  S(x_1,x_2,\ldots,x_d)\defeq \exists y_1,y_2,\ldots y_{d-1}\begin{cases}
    0&=P(y_1,y_2,\ldots, y_{d-1}),\\
    x_1&=(1-y_1)(1-y_2)\cdots(1-y_{d-1}),\\
    x_2&=(2-y_1)(2-y_2)\cdots(2-y_{d-1}),\\
    & \vdots\\
    x_d&=(d-y_1)(d-y_2)\cdots(d-y_{d-1}).
  \end{cases}
\end{align*}
A point $x:=(x_1,\ldots,x_d)\in\rel^d$ is in the set $S$ if and only
if one can find real numbers $y_1,\ldots,y_{d-1}$ for which the above
equations hold. The idea behind this definition comes from the fact
that,
for $x, y_1, \ldots, y_{d-1}$ as above,
one can construct a $d\times d$ matrix $M$ with rational entries such that
\begin{align*}
  \left(M^nx\right)_1=(n-y_1)(n-y_2)\cdots (n-y_{d-1}) \, ,
\end{align*}
for all $n\in\nat$, where $(\cdot )_1$ refers to  the first entry of the vector. Admitting the existence of such a matrix,
let $H$ be the hyperplane having normal vector $(1,0,\ldots, 0)$ and going through the origin. Then clearly $\mathcal O(M,x)$ enters $H$ at least $d-1$ times if and only if the reals $y_1,\ldots,y_{d-1}$ are distinct natural numbers, because only then is the first entry of $M^nx$---the polynomial $(n-y_1)\cdots (n-y_{d-1})$---equal to zero.

The existence of the matrix $M$ rests on the fact that the expression
$u_n = (n-y_1)\cdots (n-y_{d-1})$ (for fixed $y_1,\ldots,y_{d-1}$) can be obtained as a linear recurrence sequence of order $d$, and in turn such a linear recurrence sequence can be represented as the sequence of fixed-position entries of increasing powers of a fixed $d\times d$ matrix $M$. In the standard construction of this matricial representation, one must in addition set
$x_1=u_1$, $x_2=u_2$, \ldots, $x_d=u_d$, which is achieved through the definition of our initial semialgebraic set $S$.



It is worth noting that Theorem~\ref{th:undecidable} holds even if $k$ is fixed, due to the fact that Hilbert's tenth problem remains undecidable for a fixed number of variables. Furthermore, if $k$ is fixed to be $1$, then the problem becomes decidable in low dimensions, however even in the case where the ambient space has dimension $2$ and $k=2$, the problem does not seem to be trivially decidable.

\section{Research Directions and Open Problems}
\label{sec:future}

We have presented an overview of the state of the art regarding
decidability and solvability of a range of algorithmic problems for
discrete linear dynamical systems, focussing on reachability,
model-checking, and invariant-generation questions. In the case of model
checking in particular, we have painted an essentially complete
picture of what is achievable, both unconditionally as well as
relative to Skolem oracles. We pointed out that extending the existing
results further runs up against formidable mathematical obstacles
(longstanding open problems in number theory); the work presented here
therefore appears to lie at the very frontier of what is attainable,
barring major breakthroughs in mathematics.

A central open question is whether Skolem oracles can actually be
devised and implemented. As remarked earlier, a certifying algorithm for the Simple-Skolem
Problem has recently been proposed, with termination relying on
classical (and widely believed) number-theoretic conjectures. Whether similar
algorithms can be obtained for the general Skolem Problem is the
subject of active research. 

An even more difficult question is whether oracles for the Positivity
Problem (or Simple-Positivity Problem) can be devised. This would
enable one to circumvent the mathematical obstacles mentioned earlier,
and allow one to substantially extend the scope of semialgebraic model
checking for linear dynamical systems. For the time being, this goal
appears to be well out of reach.

In the present paper we have entirely confined ourselves to matters of
decidability, largely on account of space limitations, but also
because the complexity-theoretic picture is not nearly as clear-cut as
its decidability counterpart. Providing a comprehensive account of the complexity
landscape for linear dynamical systems would be an interesting and
promising research direction.

\appendix

  \section{Prefix-Independent Model Checking for LDS}
\label{appA}

The goal of this appendix is to exhibit boundaries on the extent to
which Theorem~\ref{thm:mc3} can be improved. More precisely, we show
that the ability to solve the model-checking problem for arbitrary LDS against
prefix-independent MSO specifications making use of semialgebraic
predicates in ambient space $\R^4$ would necessarily entail major
breakthroughs in Diophantine approximation.

We build upon the framework developed
in~\cite[Sec.~5]{DBLP:conf/soda/OuaknineW14}. To this end, consider the class
of order-$6$ rational LRS of the form
\[ u_n = -n + \frac{1}{2} (n - ri)\lambda^n + \frac{1}{2}(n +
          ri)\overline{\lambda}^n = r \operatorname{Im}(\lambda^n) - n
          (1 - \operatorname{Re}(\lambda^n)) \, ,
	\]
where $\lambda \in \mathbb{Q}(i)$ and $|\lambda| = 1$, and $r \in
\mathbb{Q}$. Let us write $\mathcal{L}$ to denote this class of LRS\@.

It is shown in~\cite{DBLP:conf/soda/OuaknineW14} that solving the
\emph{Ultimate Positivity Problem} for LRS in $\mathcal{L}$, i.e.,
providing an algorithm which, given an LRS $\langle u_n
\rangle_{n=0}^\infty \in \mathcal{L}$, determines whether there exists
some integer $N$ such that, for all $n \geq N$, $u_n \geq 0$, would
necessarily entail major breakthroughs in the field of Diophantine
approximation. The purpose of the present section is to reduce the
Ultimate Positivity Problem for LRS in $\mathcal{L}$ to the 
prefix-independent semialgebraic MSO model-checking problem for $4$-dimensional LDS\@.

Given $\lambda$ and $r$ as above, let \[
		M = \begin{bmatrix}
			\operatorname{Re}(\lambda) & -\operatorname{Im}(\lambda) & 1 & 0 \\
			\operatorname{Im}(\lambda) & \operatorname{Re}(\lambda) & 0 & 1 \\
			0 & 0 & \operatorname{Re}(\lambda) & -\operatorname{Im}(\lambda) \\
			0 & 0 & \operatorname{Im}(\lambda) & \operatorname{Re}(\lambda) \\
		\end{bmatrix} \text{ and }
		x = \begin{bmatrix}
			1\\1\\1\\1
		\end{bmatrix}.
		\] 	
		Observe that $M$ has rational entries. We have that $$M^n x = \begin{bmatrix}
			\operatorname{Re}(\lambda^n) - \operatorname{Im}(\lambda^n) + n\operatorname{Re}(\lambda^{n-1}) - n\operatorname{Im}(\lambda^{n-1})\\
			\operatorname{Im}(\lambda^n) + \operatorname{Re}(\lambda^n) + n\operatorname{Im}(\lambda^{n-1})  + n\operatorname{Re}(\lambda^{n-1}) \\
			\operatorname{Re}(\lambda^n) - \operatorname{Im}(\lambda^n) \\
			\operatorname{Im}(\lambda^n) + \operatorname{Re}(\lambda^n)	
		\end{bmatrix}.$$
		As a semialgebraic target consider the set $S =  \{x: p(x) > 0\}$, where
		\[
		p(x_1, x_2, x_3, x_4) = \frac{r}{2}(x_4 - x_3) - \frac{x_1-x_3}{\operatorname{Re}(\lambda^{-1})x_3 - \operatorname{Im}(\lambda^{-1})x_4}\left(1 - \frac{x_3 + x_4}{2}\right).
		\]
                We now have that
                $p(M^n x) = r \operatorname{Im}(\lambda^n) - n (1 -
                \operatorname{Re}(\lambda^n))$, and that $\langle u_n
\rangle_{n=0}^\infty$ is
                ultimately positive if and only if the orbit of $x$
                under $M$ eventually gets trapped in $S$. This can be
                expressed by the prefix-independent LTL formula
                $\varphi = \mathbf{F} \, \mathbf{G} \, S$.

\bibliographystyle{splncs04}
\bibliography{references2}

\end{document}

%% file: macros.tex
\newcommand{\R}{\mathbb{R}}
\newcommand{\Q}{\mathbb{Q}}
\newcommand{\N}{\mathbb{N}}

\newcommand{\C}{\mathbb{C}}

\newcommand{\T}{\mathcal{T}}
\newcommand{\PR}{\mathcal{P}}
\newcommand{\SA}{\mathcal{S}}
\newcommand{\CO}{\mathcal{C}}

\newcommand{\set}[1]{\{#1\}}

\newcommand{\nat}{\mathbb{N}}
\newcommand{\intg}{\mathbb{Z}}
\newcommand{\rel}{\mathbb{R}}

\newcommand\defeq{\mathrel{\overset{\makebox[0pt]{\mbox{\normalfont\tiny\sffamily def}}}{=}}}

\newcommand{\aut}{\mathcal{A}}
